\documentclass[11pt]{amsart}
\usepackage{amssymb,amsmath}
 \newcommand{\Q}{{\mathbb Q}}
\newcommand{\Qbar}{{\overline{\Q}}} \newcommand{\Fbar}{{\overline{F}}}
\newcommand{\Z}{{\mathbb Z}} 
 
\newcommand{\Proj}{{\mathbf P}} 
\newcommand{\TT}{{\mathcal T}} 
 \newcommand{\EE}{{\mathcal E}}
\newcommand{\LL}{{\mathcal L}} \newcommand{\C}{{\mathbb C}}

 \newcommand{\pp}{{\mathfrak p}}
\newcommand{\PP}{{\mathfrak T}} 
\newcommand{\QQ}{{\mathfrak Q}} 
 
 \newcommand{\F}{{\mathbb F}}

\newcommand{\Hom}{\operatorname{Hom}}
\newcommand{\rank}{\operatorname{rank}}
\newcommand{\Aut}{\operatorname{Aut}}
\newcommand{\Gal}{\operatorname{Gal}}
\newcommand{\Jac}{\operatorname{Jac}}
\newcommand{\Frac}{\operatorname{Frac}}

\newcommand{\isom}{\cong} \newcommand{\tensor}{\otimes}
 
\newcommand{\ord}{{\operatorname{ord}}}

\newcommand{\comment}[1]{} 
\newtheorem{thm}{Theorem}[section]
\newtheorem{lem}[thm]{Lemma}

\newtheorem{prop}[thm]{Proposition}

\theoremstyle{definition}
 
\newtheorem{defn}{Definition}
\newtheorem{rem}{Remark}
\begin{document}
\title[Hilbert's Tenth Problem for $p$-adic function fields]{Hilbert's
  Tenth Problem for function fields of varieties over number fields
  and $p$-adic fields} \thanks{As she was completing this
  paper, the author learned that Laurent Moret-Bailly had
  independently obtained the same result.}

\author{Kirsten Eisentr\"ager} \address{Department of Mathematics,
  University of Michigan, 530 Church Street, Ann Arbor, MI 48109}
\email{eisentra@umich.edu}

\begin{abstract}
  Let $k$ be a subfield of a $p$-adic field of odd residue
  characteristic, and let $\LL$ be the function field of a
  variety of dimension $n\geq 1$ over $k$. Then Hilbert's Tenth Problem
  for $\LL$ is undecidable. In particular, Hilbert's
  Tenth Problem for function fields of varieties over number fields of
  dimension $\geq1$ is undecidable.
\end{abstract}

\maketitle 

\section{Introduction}\label{intro}
Hilbert's Tenth Problem in its original form was to find an algorithm
to decide, given a polynomial equation $f(x_1,\dots,x_n)=0$ with
coefficients in the ring $\Z$ of integers, whether it has a solution
with $x_1,\dots,x_n \in \Z$.  
Matijasevi{\v{c}} (\cite{Mat70}), building on earlier work by Davis,
Putnam, and Robinson (\cite{DPR61}), proved that no such algorithm
exists, {\it i.e.}\ Hilbert's Tenth Problem is undecidable.

Since then, analogues of this problem have been studied by asking the
same question for polynomial equations with coefficients and solutions
in other commutative rings $R$.  We refer to this as {\em Hilbert's
  Tenth Problem over $R$}.  Perhaps the most important unsolved
question in this area is Hilbert's Tenth Problem over the field of
rational numbers. The function field analogue, namely Hilbert's Tenth
Problem for the function field $k$ of a curve over a finite field, is
undecidable.  This was proved by Pheidas for $k=\F_q(t)$ with $q$
odd~(\cite{Phei91}), and by Videla~(\cite{Vi94}) for $\F_q(t)$ with
$q$ even.  Shlapentokh~(\cite{Sh2000}) generalized Pheidas' result to
finite extensions of $\F_q(t)$ with $q$ odd and to certain function
fields over possibly infinite constant fields of odd characteristic,
and the remaining cases in characteristic $2$ are treated in
\cite{Eis2002}.  Hilbert's Tenth Problem is also known to be
undecidable for several rational function fields of characteristic
zero: In 1978 Denef proved the undecidability of Hilbert's Tenth
Problem for rational function fields over formally real
fields~(\cite{Den78}), and he was the first to use rank one elliptic
curves to prove undecidability. Kim and Roush~(\cite{KR92}) showed
that the problem is undecidable for the purely transcendental function
field $\C(t_1,t_2)$ and in \cite{Eis} their approach was generalized
to finite extensions of $\C(t_1,\dots,t_n)$ for $n\geq 2$.  Kim and
Roush~(\cite{KR95}) proved that the problem was undecidable for
rational function fields $k(t)$, where $k$ is a subfield of a $p$-adic
field of odd residue characteristic. In this paper we will generalize
their result to finite extensions of the rational function field in
$n$ variables over $k$ with $n \geq 1$. In particular, we show that
Hilbert's Tenth Problem for function fields of varieties over number
fields of dimension $\geq 1$ is undecidable.

In Hilbert's Tenth Problem the coefficients of the equations have to
be input into a Turing machine, so we restrict the coefficients to a
subring $A$ of $R$ which is finitely generated as a $\Z$-algebra.  We
say that {\em Hilbert's Tenth Problem for $R$ with coefficients in
  $A$} is undecidable if there is no algorithm that decides whether or
not multivariate polynomial equations with coefficients in $A$ have a
solution in $R$.
\\
Our theorem considers fields which are extensions of the rational
function field $\Q_p(\tau)$. Since $\Q_p(\tau)$ is uncountable, its
elements cannot be coded into a Turing machine. So just to get a
nontrivial problem, we have to restrict the ring of coefficients as
explained above. Let $k$ be a subfield of a $p$-adic field and let
$\LL$ be a finite extension of the rational function field
$k(\tau,\tau_2,\dots,\tau_n)$, which is given via the minimal polynomial of
a generator $\alpha$ over $k(\tau, \tau_2, \dots, \tau_n)$. (For simplicity of
notation, we assume that $\LL / k(\tau, \tau_2, \dots, \tau_n)$ is given to
us in terms of {\em one} generator $\alpha$.) We will choose the ring of
coefficients in terms of the given transcendentals $\tau,\tau_2, \dots,
\tau_n$ and $\alpha$, and we want to choose this ring as small as possible.
We will define a field $\kappa$ such that $\kappa(\tau)$ contains the
coefficients of the minimal polynomial of $\alpha$, and we will choose the
ring of coefficients to be a subring of $\kappa(\tau)$.
The field $\kappa$ will be defined in Section~\ref{definekappa}. We will
prove the following theorem:
\begin{thm}\label{Theorem}
  Let $k$ be a subfield of a finite extension of $\Q_p$ with $p$ odd.
Let $\LL$ be a finite extension of the rational function field
$k(\tau,\tau_2, \dots,\tau_n)$.
There exists a finite set $\{c_1,\dots,c_{\ell}\}$ of
  elements of $\kappa(\tau)$, not all constant, such that Hilbert's Tenth
  Problem for $\LL$ with coefficients in $\Z[c_1,\dots,c_{\ell}]$ is
  undecidable.
\end{thm}

 {\bf Notation:} In the following we will let $A_0$ be
the ring of coefficients of Theorem~\ref{Theorem}, and
$k(\tau_1,\dots,\tau_n)$ will denote the field of rational functions over
$k$ in $n$ variables $\tau_1, \dots, \tau_n$. We refer to a subfield $k$ of
a finite extension $F$ of $\Q_p$ as a {\em $p$-adic field}, and we
assume that $k$ is given together with an embedding into $F$. The
$p$-adic valuation on $F$ induces a valuation on $k$, which we
normalize so that the value group of $k$ is $\Z$.
For an integral domain $R$, we denote its field of fractions
of $R$ by $\Frac(R)$. 
\subsection{Idea of proof}
First we will define two notions that will appear frequently in
the remainder of this paper.
\begin{defn}
  1. If $R$ is a commutative ring, a {\em diophantine equation over
    $R$} is an equation $f(x_1,\dots,x_n)=0$ where $f$ is a polynomial
  in the variables $x_1, \dots, x_n$ with coefficients in
  $R$.\\
  2. A subset $S$ of $R^k$ is {\em diophantine over R} if there exists a
  polynomial\linebreak $f(x_1,\dots,x_k,$ $y_1,\dots,y_m) \in
  R[x_1,\dots,x_k,y_1,\dots,y_m]$ such that \[S=\{(x_1, \dots, x_k)
  \in R^k: \exists y_1, \dots,y_m \in R,\;
  (f(x_1,\dots,x_k,y_1,\dots,y_m)=0)\}.\]

  Let $A$ be a subring of $R$ and suppose that $f$ can be chosen such
that its coefficients are in $A$. Then we say that $S$ is {\em diophantine
over $R$ with coefficients in $A$}.

\end{defn}

We will prove Theorem~\ref{Theorem} by constructing a diophantine
model of the integers with addition and multiplication over $\LL$. A
{\em diophantine model} is defined as follows:
\begin{defn} A {\em diophantine model} of $\langle\Z,0,1;+,\cdot
  \rangle$ over $\LL$ is a diophantine subset $S \subseteq \LL^m$
  equipped with a bijection $\phi: \Z \to S$ such that under $\phi$,
  the graphs of addition and multiplication correspond to diophantine
  subsets of $S^3$.

  Let $A$ be a subring of $\LL$. A {\em diophantine model} of
  $\langle\Z,0,1;+,\cdot \rangle$ over $\LL$ {\em with coefficients in
    $A$} is a diophantine model of $\langle\Z,0,1;+,\cdot \rangle$,
  where in addition $S$ and the graphs of addition and multiplication
  are diophantine over $\LL$ with coefficients in $A$.
\end{defn}

Since Hilbert's Tenth Problem over $\Z$ is undecidable, it follows
that the structure $\langle\Z,0,1;+,\cdot \rangle$ has an undecidable
existential theory.  Hence constructing a diophantine model of
$\langle\Z,0,1;+,\cdot \rangle$ over $\LL$ with coefficients in
$A_0=\Z[c_1,\dots,c_{\ell}]$ is enough to prove that Hilbert's Tenth
Problem for $\LL$ with coefficients in $A_0$ is undecidable. We have
to check that the diophantine definition of the set $S$ which is in
bijection to $\Z$ and the diophantine definitions of addition and
multiplication have coefficients in $A_0$.  We specify the ring $A_0$
in Sections~\ref{coefficients} and \ref{proof}.

We will use the rational points on a rank one elliptic curve over
$\LL$ as our set $S$. This elliptic curve is constructed in
Section~\ref{ellipticsection}.  In Section~\ref{density} we will
generalize a theorem in \cite{KR95} to construct a diophantine set over
$\LL$ whose intersection with $\Q$ is dense in any finite product of
$p$-adic fields. In Section~\ref{quadraticforms} we prove a result
about quadratic forms that will be needed in the proof of
Theorem~\ref{Theorem}. In Section~\ref{coefficients} we address the
ring of coefficients $A_0$, and in Section~\ref{proof} we prove
Theorem~\ref{Theorem}. We will first prove Theorem~\ref{Theorem} when
$\LL/k$ has transcendence degree one, and then generalize it to higher
transcendence degree.

{\bf Note:} When she was completing the proof of
Theorem~\ref{Theorem} the author worked with an earlier version of
\cite{MB01} that did not contain the section on $p$-adic fields.

\section{Preliminaries}
We need two general facts about diophantine equations that allow us to
combine several diophantine equations into one.

\begin{lem}\label{combine}
  Let $R$ be an integral domain. Let $A$ be a subring of $R$, and
  assume that $\Frac{R}$ does not contain the algebraic closure of $\Frac(A)$.
  Then for each system $f_1(x_1, \dots, x_n)=0, \dots,
  f_k(x_1,\dots,x_n)=0$ of diophantine equations with coefficients in
  $A$ there exists a single diophantine equation $g(x_1,\dots, x_n)=0$
  with coefficients in $A$ such that the system of the $f_i$'s has a
  solution in $R$ if and only if $g$ has a solution in $R$.
\end{lem}
\begin{proof}
  We will show how to combine two equations into one, which is enough.
  Let $h(x)$ be a polynomial in one variable with coefficients in $A$
  which has no zero in $R$.  Let $\tilde{h}(x,y)$ be the
  homogenization of $h$.  Then for all $x$ and $y$ in
  $R$, $\tilde{h}(x,y)=0$ if and only if $x=0$ and $y=0$. Hence for
  $\vec{x} \in R^n$
  \[ \left( f_1(\vec{x})=0 \land f_2(\vec{x})=0 \right) \iff
  \left(\tilde{h}(f_1(\vec{x}),f_2(\vec{x}))=0\right).
\]
\end{proof}
\begin{rem}
{\em Similarly, if $R$ is an integral domain and $f_1=0$, and $f_2=0$ are
diophantine equations with coefficients in some subring $A$ of $R$, then
\[
f_1=0 \lor f_2 =0 \iff f_1\cdot f_2 =0,
\]
and $f_1\cdot f_2$ has coefficients in $A$.}
\end{rem}
\section{Algebraic function fields}\label{definekappa}
An {\em algebraic function field} in one variable over $F$ is a field
$K$ containing $F$ and at least one transcendental element $\tau$ such
that $K/F(\tau)$ is a finite algebraic extension, and such that $F$ is
algebraically closed in $K$. The field $F$ is the {\em constant field}
of $K$.
Whenever $K / F$ is an algebraic function field, we fix an algebraic
closure $\overline{K}$ of $K$.  For any field $E \subseteq \overline{K}$, we set $K E$
equal to the {\em compositum} of $K$ and $E$ inside $\overline{K}$.

We first need a general theorem about extensions of function fields.
\begin{thm}\label{resextension}
  Let $K/F$ be a function field of characteristic zero with constant
  field $F$. Let $E$ be an extension of $F$, and let $L=KE$. Let
  $\mathcal{T}$ be a prime of $L$ lying above a prime $\mathfrak{T}$
  of $K$. Let $L_{\mathcal{T}}$ and $K_{\mathfrak{T}}$ be the
  corresponding residue fields.
  \\
  (1) If $E/F$ is finite, then $L_{\mathcal{T}}$ is the composite of
  the two subfields $K_{\mathfrak{T}}$ and $E$.
  \\
  (2) If $E$ is algebraically closed in $L$ and $E \cap K =F$, then
  $L_{\mathcal{T}} =K_{\mathfrak{T}} E$.

\end{thm}
\begin{proof}
The first part is proved in \cite[p.\ 106]{Rosen}, and the second part
is proved in \cite[p.\ 128]{Deuring}.
\end{proof}
\subsection{Definition of the field $\kappa$}\label{kappa}
Assume that $k$ is a field of characteristic zero and and that $\LL/
k$ is an algebraic function field with constant field $k$. We will
assume that $\LL$ is specified as $k(\tau)(\alpha)$, where $\tau$ is
transcendental over $k$ and $\alpha$ generates $\LL$ over $k(\tau)$. Let
$\beta_1, \dots, \beta_n \in k(\tau)$ be the coefficients of the minimal
polynomial of $\alpha$.  Then $\beta_i = p_i(\tau)/q_i(\tau)$, with $p_i,q_i \in
k[\tau]$. Let $\kappa$ be the subfield of $k$ generated by the coefficients
of all the $p_i,q_i, i=1, \dots, n$.  Then $\kappa$ is a finitely
generated extension of $\Q$, and $\kappa(\tau)$ contains the coefficients of
the minimal polynomial of $\alpha$.  Let $K$ be the subfield of $\LL$
defined by $K:= \kappa (\tau,\alpha)$.  By construction, the field $K$ is an
algebraic function field with constant field $\kappa$.
\begin{prop}\label{ff1}
  The field $\LL/k$ is a constant field extension of $K/ \kappa$, {\it
    i.e.}\ $K k = \LL$, and $k \cap K = \kappa$.
\end{prop}
\begin{proof} We have $Kk = \LL$ by construction. It remains to show
  that $k \cap K = \kappa$.  We will show this by showing that $k$ is
  linearly disjoint from $K$ over $\kappa$. By \cite[Lemma 3, p.\
  123]{Deuring} applied to $\kappa \subset \kappa(\tau) \subset K$ and
  $\kappa \subset k$, it suffices to show that $\kappa(\tau)$ is
  linearly disjoint from $k$ over $\kappa$, and that
  $k(\kappa(\tau))=k(\tau)$ is linearly disjoint from $K$ over
  $\kappa(\tau)$. Since $\tau$ is transcendental over $k$,
  $\kappa(\tau)$ is linearly disjoint from $k$ over $\kappa$
  (\cite[Lemma 2 (a), p.\ 122]{Deuring}).  By construction of $\kappa$
  and $K$, $\LL= k(\tau)(\alpha)$, $K=\kappa(\tau)(\alpha)$, and
  $[k(\tau)(\alpha):k(\tau)]=[\kappa(\tau)(\alpha):\kappa(\tau)]$. Hence
  $k(\tau)$ is linearly disjoint from $K$ over $\kappa(\tau)$ by
  \cite[Lemma 1, p. 122]{Deuring}.
\end{proof}
\begin{prop}\label{ff2}
  Let $K/ \kappa(\tau)$ and $\LL / k(\tau)$ be as before, and assume
  that there is a prime $\TT$ of $K$ above $\tau$ which is
  unramified. Then there exists a prime $\TT'$ of $\LL$ above $\tau$
  which is unramified.
  \\
  Moreover, there exists a finite extension $\kappa_1$ of $\kappa$
  such that in the compositum of $\LL$ and $k\kappa_1$, the residue
  field of any prime extending $\TT'$ is $k \kappa_1$.
\end{prop}
\begin{proof}
  By Proposition~\ref{ff1}, the extension $\LL / k$ is a constant
  field extension of $K / \kappa$. Hence the prime $\TT'$ of $\LL$
  extending $\TT$ is unramified (\cite[p.\ 113]{Deuring}), and hence
  $\TT'$ over $\tau$ is unramified.
  \\
  Let $K_{\TT}$ be the residue field of the prime $\TT$ of $K$. Then
  $K_{\TT}$ is a finite extension of $\kappa$. By
  Theorem~\ref{resextension}(2), the residue field $\LL_{\TT'}$ of the
  prime $\TT'$ above $\TT$ is $K_{\TT}k$.  Similarly, let $k'$ be a
  finite extension of $k$, and let $\mathfrak{Q}$ be a prime of $\LL
  k'$ extending $\TT'$. By Theorem~\ref{resextension}(1), the residue
  field of $\mathfrak{Q}$ is the compositum of $\LL_{\TT'}$ and $k'$.
  Now let $\kappa _1$ be a finite (normal) extension of $\kappa$, such
  that $K_{\TT} \kappa_1 = \kappa _1$. Then by the above arguments,
  $\kappa_1$ has the right properties.
\end{proof}

\subsection{Higher transcendence degree}\label{highertrdeg}
When $\LL$ is a finite extension of\linebreak $k(\tau,\tau_2
\dots,\tau_n)$ (with $k$ algebraically closed in $\LL$), which is
given as $\LL= k(\tau, \tau_2, \dots, \tau_n)(\alpha)$, then the
coefficients $\beta_i$ of the minimal polynomial of $\alpha$ over
$k(\tau, \tau_2, \dots, \tau_n)$ are elements of $k(\tau, \tau_2,
\dots, \tau_n)$. So each $\beta_i$ is of the form \[\beta_i =
p_i(\tau, \tau_2, \dots, \tau_n)/q_i(\tau, \tau_2, \dots, \tau_n)
\text{ with } p_i,q_i \in k[\tau, \tau_2, \dots, \tau_n].\] Let
$\kappa_0$ be the subfield of $k$ which is generated by the
coefficients of the $p_i,q_i$.

\begin{prop}
Let $k_1$ be the algebraic closure of $k(\tau_2, \dots, \tau_n)$ in $\LL$.
There exists a finite extension $\kappa$ of $\kappa _0(\tau_2,\dots,\tau_n)$ and
a finite extension $K$ of $\kappa(\tau)$ such that the algebraic function
field $\LL / k_1$ is a constant field extension of $K / \kappa$.
\end{prop}
\begin{proof}
The field $k_1$ is a finite extension of $k(\tau_2, \dots, \tau_n)$.
The coefficients of the minimal polynomial of $\alpha$ over $k_1(\tau)$
are algebraic over $\kappa_0(\tau ,\tau_2, \dots, \tau_n)$, and generate some
finite extension $K _1$ of $\kappa_0(\tau ,\tau_2, \dots, \tau_n)$ which is
contained in $k_1(\tau)$. Let $K:=
K_1(\alpha) \subseteq \LL$, and let $\kappa$ be the algebraic closure of $\kappa_0(\tau_2,
\dots, \tau_n)$ in $K$.
\\
Then with the same argument as in the proof of Proposition~\ref{ff1},
the algebraic function field $\LL / k$ is a constant field extension
of $K / \kappa$.
\end{proof}
\begin{rem}\label{higher}
  Exactly the same proof as the proof of Proposition~\ref{ff2} shows
  that Proposition~\ref{ff2} also holds for the extension $\LL /
  k_1(\tau)$ of $K / \kappa(\tau)$.
\end{rem}
\section{Elliptic curve setup}\label{ellipticsection}
To construct a diophantine model of $\langle\Z,0,1;+,\cdot \rangle$ over $\LL$ with
coefficients in $A_0$ we need a diophantine set $\mathcal{S}$ and a
bijection $\Z\to \mathcal{S}$. We will choose as our set $\mathcal{S}$
the $\LL$-rational points on an elliptic curve $\EE_0$, and so we need
an elliptic curve $\EE_0$ over $\LL$ of rank one. The following
theorem uses a theorem by Moret-Bailly
(\cite[Theorem 1.8]{MB01}) and allows us to construct elliptic
curves of rank one:
\begin{thm}\label{elliptic}
  Let $k$ be a field of characteristic zero, and let $\LL$ be a finite
  extension of the rational function field $k(\tau)$. Let $\kappa$ and
  $K$ be as in Proposition~\ref{ff2}, and let $E$ be an elliptic curve
  over $\Q$
  without complex multiplication and with Weierstrass equation
  $$y^2=x^3+ax+b,$$ where $a,b \in \Q$, $b \neq 0$. Then there exists a non-constant
  element $T \in K$ such that the elliptic curve given by the affine
  equation $$\EE: (T^3 + aT +b)Y^2=X^3+aX+b$$ has rank one over $\LL$
  with generator $(T,1)$ modulo 2-torsion. Moreover, $T$ can be chosen
  such that the extension $\LL / k(T)$ is unramified above the primes
  $T,T^{-1}$ of $k(T)$.
\end{thm}
\begin{proof}
Let $C_0$ be a smooth projective geometrically connected curve over
$\kappa$ with function field $K$.
\\
To get the desired element $T$, pick an ``admissible element'' $f \in
K$ (\cite[Definition 1.5.2]{MB01}), pick an element $\lambda \in \Z\;\cap$
GOOD$(\kappa)$\cite[Theorem 1.8]{MB01}, and let $T:= \lambda\cdot f$. Then
$\EE(K)$ is generated by $(T,1)$.  Since $T$ is admissible in the
sense of \cite[Definition 1.5.2]{MB01}, it follows that $T:C_0 \to
\mathbb{P}^{1}_{\kappa}$ is \'etale above $0$ and $\infty$.
\\
Moreover, the group $\EE(\LL)$ is generated by $(T,1)$: Indeed, the
field $k$ is an extension of $\kappa$, and by \cite[Corollary
1.5.5(ii)]{MB01}, GOOD$(k) \cap \kappa =$ GOOD$(\kappa)$, so $T\in$ GOOD$(k)$. By the
definition of ``GOOD'', this means that the natural inclusion
$\EE(k(T)) \hookrightarrow \EE(\LL)$ is a bijection, so $\EE(\LL)$
is generated by $(T,1)$.
\end{proof}

{\bf Note:} Our notation follows Moret-Bailly's equivalent setup in his
preprint of \cite{MB01} from December 2003: We assume that the
polynomial $R(t)$ defining $\Gamma$ in 1.4.4. is without multiple roots
and satisfies $R(0) \neq 0$.  We are also in the situation $\Gamma = E$, but
the double cover $\pi$ is given by the $x$-coordinate. With this
notation, we have $R(t)=P(t)$ and the twisted curve
$y^2=R(t)P(x)$ in \cite[1.4.6]{MB01} is isomorphic to $R(t)y^2 =P(x)$
(which is the twist that we use) via $(X,Y) \mapsto (X,Y/R(t))$.

{\bf Notation:} For $k,\LL,T,\EE$ as above let $P:=(T,1)$. Let
$P_m:=m \cdot (T,1)=(X_m,Y_m)$, and for $m \neq 0$ let $\psi_m:=X_m/TY_m$.

Since $\psi_m \in k(T)$, we can interpret $\psi_m$ as a function on the
projective line. We will need a proposition by Denef, which determines
$\psi_m(\infty)$.
\begin{prop}\label{Denef}
  The function $\psi_m$ takes the value $m$ at $\infty$. {\it I.e.}, when we
  expand $\psi_m$ as a power series in $T^{-1}$, the constant term is
  $m$.
\end{prop}
\begin{proof}
This is Lemma 3.2 in \cite[p.~396]{Den78}.
\end{proof}
In the proof of
Theorem~\ref{Theorem} we will need more properties of the points on
the elliptic curve, so we will work with a specific curve that has
these properties. From now on we will fix $E_0$ to be the
smooth projective model of $$y^2 =x^3 + x+ 1.$$ Then $E_0$ has
no complex multiplication, and the point $(0,1) \in E_0(k)$ has
infinite order (\cite[496A1]{Cremona}). We will fix an element $T \in
\LL$ as in Theorem~\ref{elliptic}. Let $\EE_0$ be the elliptic curve
given by
$$\EE_0: (T^3+T+1)Y^2=X^3+X+1.$$ By our choice of $T$, a generator
for $\EE_0(\LL)$ (modulo 2-torsion) is $(T,1)$.
\begin{lem}\label{orders}
  Given $E_0,\EE_0$, let $\psi_m$ be defined as above.  Given $m,n,r \in
  \Z-\{0,1,-1\}$, let
  \begin{gather*}u:= u_{(m,n,r)}= \psi_{m} \psi_{n}-\psi_{r} + ({1}/{2})\cdot
    T^{-1}, \\ v:= v_{(m,n,r)}=
    \psi_{m} \psi_{n}-\psi_{r} + ({1}/{3})\cdot T^{-1}.
\end{gather*}
Let $\ord_T,\ord_{T^{-1}}$ be the discrete valuations on $k(T)$ associated to
$T,T^{-1}$, normalized so that the value group is $\Z$. Then
\begin{enumerate}
\item \label{stat1} $\ord_{T }(u)=-2$, $\ord_T(v) =-2$.
\item $n \cdot m=r$ if and only if $(\ord_{T^{-1}} (u) =1 \text{ or
  }\ord_{T^{-1}} (v) =1)$.\\
$ n \cdot m \neq r$ if and only if $(\ord_{T^{-1}} (u)
  = 0 \text{ and } \ord_{T^{-1}} (v) =0)$.
\end{enumerate}

\end{lem}
\begin{proof}$\left.\right.$
  (1) When we reduce the equation of the curve $\EE_0$ modulo $T$ we
  just obtain $E_0$, so the reduction of $\EE_0/k(T)$ modulo $T$ gives us
  the nonsingular curve $E_0/k$. We have a map $\pi: \EE_0 \to E_0$ that
  sends a point $Q \in E_0(k(T))$ to a point $\tilde Q$, its
  reduction modulo $T$, and this map is a group homomorphism. The
  reduction of the point $P=(T,1)$ on $\EE_0$ is the point $(0,1)$ on
  $E_0$, and since $(0,1)$ has infinite order, this means that no
  non-zero multiple of $P$ can map to $\mathbf{O}$. Hence
  $X_{m},Y_{m}$ have nonnegative order at $T$ for all $m \in
  \Z-\{0\}$. Since the reduction of $P$ has infinite order it follows
  that $Y_m$ has order zero at $T$.  If $X_m$ has positive order at
  $T$, then under $\pi$ it gets mapped to a point on $E_0$ whose
  $x$-coordinate is zero. The only such points on $E_0$ are $(0,1)$ and
  $(0,-1)$ which are the images of $P_1$ and $P_{-1}$ respectively.
  Since $(0,1)\in E_0(k)$ has infinite order and since $\pi$ is a group
  homomorphism this implies that no other multiples of $P$ can reduce
  to $(0,\pm 1)$. So for $m \in \Z-\{0,1,-1\}$ we have that $X_m,Y_m$
  have order 0 at $T$.  Hence for all $m \in \Z-\{0,1,-1\}$, $\psi_m$
  has order $-1$ at $T$, and so $u=\psi_{m} \psi_{n}-\psi_{r} +
  (1/2) \cdot T^{-1}$ has order $-2$ at $T$. Similarly, $\ord_T(v)=-2$.

  (2) If $n\cdot m = r$, then by Proposition~\ref{Denef}, $\psi_{m}
  \psi_{n}-\psi_{r}$ has nonnegative order in $T^{-1}$ and the constant
  coefficient cancels, so the order at $T^{-1}$ is positive. Hence at
  least one of the power series expansion of $u$ and $v$ in $T^{-1}$
  has a linear term, and so $\ord_{T^{-1}}(u) =1$ or
  $\ord_{T^{-1}}(v)=1$.

 If $n \cdot m\neq r$, then by Proposition~\ref{Denef}, $\psi_{m}
  \psi_{n}-\psi_{r}$ has nonnegative order in $T^{-1}$ and the
  constant term in the power series expansion in $T^{-1}$ is nonzero.
  Hence $\ord_{T^{-1}}(u) =0$,  $\ord_{T^{-1}}(v) =0$.
\end{proof}
\section{A set, diophantine over $\LL$, which is dense in any finite
  product of $p$-adic fields}\label{density}
The goal of this section is to prove Theorem~\ref{densitytheorem},
which will be needed in the proof of Theorem~\ref{Theorem}.
\begin{prop}\label{exclusion}
  Let $k$ be a field of characteristic zero, and let $\LL$ be a finite
  extension of the rational function field $k(T)$. Assume that $k$ is
  algebraically closed in $\LL$. There exists a finite set
  $\mathcal{A}$ of elliptic curves over $\Q$ with the property that if
  $E/ \Q$ is an elliptic curve which is not $\Qbar$ isogenous to any
  of the curves in $\mathcal{A}$, then $E(k) = E(\LL)$.
\end{prop}
\begin{proof}
  Let $C$ be a smooth, projective, geometrically connected curve
  defined over $k$ whose function field is $\LL$. Let $E$ be an
  elliptic curve defined over $\Q$. A non-constant point $P \in E(\LL)$
  corresponds to a non-constant morphism $\alpha: C\to E$ defined over $k$.
  The morphism $\alpha$ induces a non-zero homomorphism $\beta:\Jac(C) \to
  \Jac(E) \isom E$ defined over $k$. We can decompose $\Jac(C)$ into
  simple factors over ${k}$.  In order to have a non-zero homomorphism
  $\beta:\Jac(C) \to \Jac(E)$ one of the simple factors $A_0$ of $\Jac(C)$
  has to be ${k}$-isogenous to $E$. So if $E$ is not $k$-isogenous to
  any of the $k$-simple factors of $\Jac(C)$, then $E(k) = E(\LL)$.

  If two elliptic curves $E_0,E_1$ defined over $\Q$ are both
  $k$-isogenous to a simple factor $A$ of $\Jac(C)$, then $E_0$ and
  $E_1$ are $k$-isogenous.  But then $E_0$ and $E_1$ must already be
  isogenous over $\Qbar$~\cite[Theorem 2.1]{Brian}. So requiring that
  an elliptic curve $E/ \Q$ not be isogenous to any of the simple
  factors of $\Jac(C)$ over $k$ excludes finitely many $\Qbar$ isogeny
  classes of elliptic curves defined over $\Q$.
\end{proof}
\begin{rem}
{\em  We can use a similar argument as above to prove
  Proposition~\ref{exclusion} when $\LL$ is a finite extension of the
  rational function field\linebreak
  $k(T_1,\dots,T_n)$ in $n$ variables with $k$ algebraically
  closed in $\LL$:
  
  Let $k_i:=k(T_1,\dots, \hat{T_i},\dots,T_n)$ for $i=1,\dots,n$.
  (Here $\hat{T_i}$ means that $T_i$ is omitted.) For $i=1,\dots, n$,
  let $K_i$ be the algebraic closure of $k_i$ in $\LL$ and let $C_i$
  be a smooth, projective, geometrically connected curve defined over
  $K_i$ whose function field is $\LL$.  Let $E/\Q$ be an elliptic
  curve. A non-constant point $P\in E(\LL)$ will have coordinates
  transcendental over some $K_i$ ($i\in \{1,\dots, n\}$), inducing a
  non-constant morphism $\alpha_i: C_i\to E$ defined over $K_i$. This gives
  a non-zero homomorphism $\beta_i:\Jac(C_i) \to E$ defined over $K_i$.
  As argued above, requiring that $E$ not be $K_i$-isogenous to any of
  the simple-factors of $\Jac(C_i)$ excludes finitely many $\Qbar$
  isogeny classes of elliptic curves over $\Q$.  Hence excluding all
  elliptic curves $E/ \Q$ which are $K_i$-isogenous to some factor of
  $\Jac(C_i)$ for some $i\in\{1, \dots, n\}$ still only excludes
  finitely many $\Qbar$ isogeny classes.}
\end{rem}
\begin{prop}\label{dense}
Let $E/\Q$ be an elliptic curve with global minimal Weierstrass equation
\[
E:y^2+a_1xy+a_3y = x^3 + a_2x^2 +a_4x+a_6.
\]
Assume that $E(\Q)$ is infinite. Let $S:=\{x/y:(x,y) \in E(\Q), y \neq
0\}$, and $U:=\{s_1/s_2: s_1,s_2 \in S, s_2 \neq 0\}$.  
\begin{enumerate}
\item Let $p$ be a prime. The $p$-adic closure of $S$ in $\Q_p$
  contains a neighborhood of the origin, and $U$ is dense in $\Q_p$.
\item Let $p_1, \dots, p_r$ be a finite set of primes. The closure of
  the set $S$ (embedded into $\Q_{p_1} \times \dots \times \Q_{p_r}$
  diagonally) contains a neighborhood of the origin, and $U$
  is dense in $\Q_{p_1} \times \dots \times \Q_{p_r}$.
\end{enumerate}
\end{prop}
\begin{proof}$\left. \right.$
  (1) Consider the curve $E$ as a curve over $\Q_p$, and let
  $\tilde{E}_{ns}(\F_p)$ be the nonsingular part of the reduction of
  $E$ modulo $p$. Let $P \mapsto \tilde{P}$ be the reduction map as in
  \cite[p.\ 173]{Silverman}. Let $E_1(\Q_p):= \{P \in
  E(\Q_p):\tilde{P} = \tilde{\mathbf{O}}\}$, and let $P_0 \in E(\Q)$
  be a point of infinite order. Some multiple of $P_0$ reduces to the
  identity, say $m\tilde P_0 = \tilde{\mathbf{O}}$.  Let
  $\widehat{E}/\Z_p$ be the formal group associated to $E$. Then
  $E_1(\Q_p) \isom \widehat{E}(p\Z _p)$ (as groups) via $(x,y) \mapsto
  -(x/y)$ (\cite[p.\ 175]{Silverman}). Hence the subgroup of $E(\Q_p)$
  generated by $mP_0$ corresponds to an infinite subgroup $G$ of the
  formal group.  Since the formal group associated to an elliptic
  curve is a one-dimensional compact $p$-adic Lie group, it follows
  that the closure of $G$ (and hence the closure of $S$) contains a
  neighborhood of the origin. \\
  Since the closure of $S$
  contains $p^n \Z_p$, it follows immediately that $U$ is dense in $\Q_p$.

(3) We can take a large enough multiple $mP_0$ of the point $P_0 \in
  E(\Q)$ of infinite order such that $mP_0$ reduces to
  the identity in the nonsingular part of the reduction of $E$
  modulo $p_i$ for  $i=1,\dots,r$. Let $R:=\Z_{p_1} \times \dots
  \times \Z_{p_r}$. The subgroup of $E(\Q)$ generated
  by $mP_0$ corresponds to an additive subgroup $M$ of $R$ via
  $$P=(x,y) \mapsto (x/y, \dots, x/y) .$$
  Let $\overline{M}$ be the
  closure of $M$ in $R$. Then $\overline{M}$ is stable under
  multiplication by $\Z$. By the strong approximation theorem
  (\cite[p.\ 67]{Cassels}) $\Z$ is dense in $R$, so it follows that
  $\overline{M}$ is stable under multiplication by elements of $R$. So
  $\overline{M}$ is an ideal of $R=\Z_{p_1} \times \dots \times
  \Z_{p_r}$. Then $\overline{M}=I_1\times \dots \times I_r$, with
  $I_i$ an ideal of $\Z_{p_i}$. By part (1) the $i$th projection of
  $\overline{M}$ contains a neighborhood of the origin, so all the
  $I_i$'s are nonzero ideals of $\Z_{p_i}$, {\it i.e.\ 
  }$I_i=p_i^{n_i}\Z_{p_i}$. Hence $\overline{M}$ contains a
  neighborhood of the origin, and $U$ is dense in $\Q_{p_1} \times
  \dots \times \Q_{p_r}$.
\end{proof}
\begin{thm}\label{help}
  Let $F$ be a number field, $E_0/F$ an elliptic curve without
  geometric complex multiplication. Let $F'$ be an extension of $F$.
  The set of $F$-isomorphism classes of elliptic curves $E/F$ which
  are $F'$-isogenous to $E_0$ is finite up to quadratic twist. {\it
    I.e.}, the set of possible $j$-invariants for $E$ is finite.
\end{thm}
\begin{proof}
  By replacing $F'$ with an extension we may assume that $F'$ is
  algebraically closed. Then $F' \supseteq \overline{F}$. If two elliptic
  curves over $F$ become isogenous over $F'$ then they are already
  isogenous over $\Fbar$ (\cite[Theorem 2.1]{Brian}), so we may assume
  $F'=\Fbar$.  Let $E/F$ be an elliptic curve as in the theorem, so
  $E_0$ and $E$ are $\overline{F}$-isogenous. Let $G_F:=
  \Gal(\Fbar/F)$. Since $E_0$ does not have geometric complex
  multiplication, $\Hom_{\Fbar}(E_0,E)$ is a free $\Z$-module of rank
  one. Thus, the natural continuous action by $G_F$ is through $G_F \to
  \Aut(\Hom_{\Fbar}(E_0,E))=\Z^{\times}=\langle \pm 1 \rangle$. That is, $E_0$ and $E$
  become isogenous over a quadratic extension $K$ of $F$.  Let $E'$ be
  the twist of $E$ by the quadratic character $\chi$ associated with
  $K/F$.
  \\
  We can show that over $F$, $E_0$ is isogenous either to $E$ or to
  $E'$: To see this, assume that $E_0$ is not isogenous to $E$ over
  $F$. Then the nontrivial $F$-automorphism of $K$, $\sigma$, acts by $-1$
  on $\Hom_{\Fbar}(E_0,E)$. Since $E$ is not isogenous to $E'$ over
  $F$, $\sigma$ also acts by $-1$ on $\Hom_{\Fbar}(E,E')$. Hence, after
  composing we see that $\sigma$ acts trivially on $\Hom_{\Fbar}(E_0,E')$,
  {\it i.e.\ }$E_0$ is $F$-isogenous to $E'$.
  \\
  But by a theorem of Shafarevich (see \cite[IX.6]{Silverman}) there
  are only finitely many $F$-isomorphism classes of elliptic curves
  defined over $F$ which are $F$-isogenous to $E_0$.
\end{proof}
Now we can prove the theorem that we will need in Section~\ref{proof}.
\begin{thm}\label{densitytheorem}
  Let $p_1, \dots, p_r$ be a finite set of primes.  Let $k$ be a
  subfield of a $p$-adic field, and let $\LL$ be a finite extension of
  $k(T)$. Assume that $k$ is algebraically closed in $\LL$. There
  exists a set $U_0 \subseteq k$ such that $U_0$ is diophantine over $\LL$
  with coefficients in $\Z$ and such that $U_0 \cap \Q$ is dense in
  $\Q_{p_1} \times \dots \times \Q_{p_r}$.
\end{thm}
\begin{proof}
  Let $\EE_{\eta}$ be an elliptic curve over $\Q(T)$, and let
  $\tilde{\EE} \to \Proj^1_{\Q}$ be an elliptic surface whose generic
  fiber is $\EE_{\eta}$. Assume that the $j$-invariant $j_{\tilde \EE}$
  of $\tilde \EE$ is non-constant, and that $\rank(\EE_{\eta}(\Q(T))) \geq
  1$. For all but finitely many $t \in \Proj^1(\Q)$, the specialization
  $\EE_t$ is an elliptic curve over $\Q$. By Silverman's
  specialization theorem (\cite{Sil83}), $\rank(\EE_{\eta}(\Q(T))) \leq
  \rank(\EE_t(\Q))$ for all but finitely many $t \in \Q$, and so
  $\EE_t$ has positive rank for all but finitely many $t \in \Q$.  We
  will now use Proposition~\ref{exclusion} and Theorem~\ref{help} to
  show that there exists a value $t \in \Q$ such that $\EE_t$ has
  positive rank, and such that $\EE_t(k)=\EE_t(\LL)$: Let $M$ be the
  set of all $t$ for which $\EE_t$ has positive rank and no geometric
  CM. Up to isomorphism over $\Qbar$ there are only a finite number of
  elliptic curves $E/ \Q$ with complex multiplication \cite[p.\ 
  340]{Silverman}, so since $j_{\tilde \EE}$ is non-constant, Silverman's
  theorem implies that $M$ is infinite and that $\{j(\EE_t):t \in M\}$
  is also infinite. If we want to ensure $\EE_t(k)=\EE_t(\LL)$, then
  by Proposition~\ref{exclusion} and Theorem~\ref{help} this excludes
  only finitely many $j$-invariants $j(\EE_t)$. Hence there is
  a $t \in M$ with the desired properties.

  Take such a $t \in \Q$ and a corresponding elliptic curve $\EE_t/
  \Q$. Let $$U_0:=\{(x/y)\cdot (y'/x'): (x,y) \in \EE_t(\LL),(x',y') \in
  \EE_t(\LL), y\cdot x' \neq 0\}.$$ Since the elliptic curve $\EE_t$ has
  coefficients in $\Q$, we can clear the denominators in its equation,
  and so $U_0$ is diophantine over $\LL$ with coefficients in $\Z$.
  Also $U_0 \subseteq k$, and by Proposition~\ref{dense}, part (2), $U_0\cap
  \Q$ is dense in $\Q_{p_1}\times \dots \times \Q_{p_r}$.
\end{proof}
\begin{rem} {\em This theorem also holds for fields $\LL$ which are
    finite extensions of $k(T_1,\dots,T_n)$ with $k$ algebraically
    closed in $\LL$: Let $\EE_{\eta}$ and $\tilde{\EE}$ be as in the
    proof of Theorem~\ref{densitytheorem}. By the remark after
    Proposition~\ref{exclusion} and Theorem~\ref{help}, to find an
    element $t$ with $\EE_t(k)=\EE_t(\LL)$ we only have to exclude
    finitely many $j$-invariants, and the proof proceeds exactly as
    before.}
\end{rem}
\section{Quadratic forms over function fields}\label{quadraticforms}
The following lemma deals with quadratic forms over $\LL$ and
generalizes Proposition~7 in \cite{KR95}. This lemma will be needed to
define multiplication on our set $\mathcal{S}$. Our notation for
quadratic forms follows~\cite{Lam}.
\begin{lem}\label{anisotropic}
  Let $k$ be a field of characteristic zero, and suppose there is a
  quadratic form $\langle 1,-a \rangle \langle 1,b \rangle$ over $k$, which is anisotropic
  over $k$. Let $\LL$ be a finite extension of $k(T)$, and let $\PP$
  be a prime above $T$ which is unramified. Assume that the residue
  field of $\LL$ at $\PP$ is $k$. Let $g \in k(T)$ be such that $\ord_T
  (g)$ is non-negative and even. Then one of the following two
  quadratic forms
\begin{gather}\nonumber
q_1=\langle T,-aT,-1,-g\rangle \langle 1,b\rangle\\
q_2=\langle T,-aT,-1,-ag\rangle \langle 1,b\rangle\label{E:1}
\end{gather}
is anisotropic over $\LL$.
\end{lem}
\begin{proof}
  Let $\ord_{\PP}:\LL^{*} \twoheadrightarrow \Z$ be the discrete
  valuation associated to $\PP$.  Since $\PP$ over $T$ is unramified,
  the element $T$ is a uniformizer for $\ord_{\PP}$.  Since $g\in k(T)$
  has even order in $T$, we may replace it by $T^{2n}g$ to ensure
  $g(0)$ is nonzero. Changing the coefficients of the quadratic forms
  by squares does not change the solvability.  Assume both forms
  represented in (\ref{E:1}) are isotropic over $\LL$. We will derive
  a contradiction from this.  Rewrite $q_1$ and $q_2$ as
\begin{gather}\label{EE:4}
T x_1^2 -T ax_2^2 + T bx_3^2 - T ab x_4^2= x_5^2 + bx_6^2+ 
gx_7^2  + bgx_8^2\\
\label{EE:5}
T y_1^2 -T ay_2^2 + T by_3^2 - T ab y_4^2= 
y_5^2 + by_6^2+ agy_7^2 + bagy_8^2
\end{gather}
  We can take a solution $(x_1, \dots, x_8)$ of $q_1$ in
$\LL$ such that $\ord_{\PP}(x_i) \geq 0$ and such that
$\ord_{\PP}(x_i)=0$ for some $i$.  Similarly we can take a solution
$(y_1, \dots, y_8)$ of $q_2$ in $\LL$ such that $\ord_{\PP}(y_i) \geq
0$ and such that $\ord_{\PP}(x_i)=0$ for some $i$.  Reduce
(\ref{EE:4}) and (\ref{EE:5}) modulo $\PP$ for these solutions.  Let
$g(0) =\ell.$ After reducing modulo $\PP$ the right-hand side of
(\ref{EE:4}) and (\ref{EE:5}) becomes $\langle 1, a^e
\ell\rangle\langle1,b\rangle$ $(e \in \{0,1\})$, which is a quadratic
form over $k$ by our assumptions on the residue field at the prime
$\PP$.  Suppose that after reducing modulo $\PP$ the right-hand side of
(\ref{EE:4}) and (\ref{EE:5}) is isotropic over $k$. The quaternion
algebras associated to $\langle1,a \ell \rangle\langle1,b \rangle$ and
$\langle1,\ell \rangle\langle1,b \rangle$ are
$\left(\frac{-b,-a\ell}{k} \right)$ and $\left(\frac{-b,-\ell}{k}
\right)$. (See Definitions \ref{quaternion} and \ref{def2} in the
appendix.) Since $\langle1,a \ell \rangle\langle1,b \rangle$ and
$\langle1,\ell \rangle\langle1,b \rangle$ are isotropic over $k$, this
implies that the quaternion algebras $\left(\frac{-b,-a\ell}{k}
\right)$ and $\left(\frac{-b,-\ell}{k} \right)$ are split over $k$
(see Proposition~\ref{split}).  But this implies that their tensor
product is isomorphic to a matrix algebra as well, and by
Proposition~\ref{multiply}, this tensor product is
$$
\left(\frac{-b,-a\ell}{k} \right)\tensor \left(\frac{-b,-\ell}{k}
\right) \isom \left(\frac{-b,a\ell^2}{k} \right)\tensor M_2(k).$$
This implies that
$\left(\frac{-b,a\ell^2}{k} \right)$ is split over $k$.  By
Proposition~\ref{split} from the appendix its associated norm form
$\langle1,b,-a\ell^2 -ab\ell^2 \rangle$ is isotropic over $k$, which
means that $\langle1,b,-a -ab \rangle$ is isotropic over $k$,
contradicting our assumptions made in the statement of the lemma.

Therefore, the right-hand side modulo $\PP$ is anisotropic for some $e
\in \{0,1\}.$ We may assume that the right-hand side of $q_1$ is
anisotropic modulo $\PP$. This can only happen if $\ord_{\PP}(x_i)>0$
for $i=5,\dots 8$. Let $\tilde{x_i} = x_i/T$ for $i=5,\dots 8$. Since
$\ord_{\PP}{T}=1$, $\ord_{\PP}(\tilde{x_i}) \geq 0$ for
$i=5,\dots,8$. We can rewrite (\ref{EE:4}) as
\[ x_1^2 -ax_2^2 + bx_3^2 - ab x_4^2= T(\tilde{x_5}
^2 + b\tilde{x_6}^2+ g\tilde{x_7}^2
  + bg\tilde{x_8}^2) 
\]
If we reduce modulo $\PP$ then we get $$\overline{x}_1^2
-a\overline{x}_2^2 + b\overline{x}_3^2 - ab \overline{x}_4^2 =0.$$
Since $\langle 1 ,-a \rangle \langle 1 , b\rangle$ was assumed to be
anisotropic over $k$ the same argument as before implies that
$\ord_{\PP}(x_i)>0$ for $i=1,\dots, 4$. That means all the $x_i$
($i=1,\dots,8$) in the solution of $q_1$ satisfy $\ord_{\PP}(x_i)>0$,
contradicting our choice of the $x_i$.
\end{proof}
\section{Enlarging the constant field and coefficients of
  equations}\label{coefficients}
%
We say that a subfield $k$ of a $p$-adic field satisfies {\em
  Hypothesis $(\mathcal{H})$} \cite[p.\ 92]{KR95}, if the
following conditions are satisfied:

There exists a four-dimensional anisotropic quadratic form $q$ over
$k$,
\[
q = \langle 1,a\rangle \langle1, \pp \rangle= x^2 + \pp y^2 + az^2+ a\pp w^2.
\] 
We require that $\pp \in k$ is an element of odd valuation, which is
algebraic over $\Q$. The element $a$ is a $2^r$-th root of unity for
some $r\geq1$, and $k$ contains a square root $i$ of $-1$. We also
require that $q$ is locally isotropic at all 2-adic primes of
$\Q(i,a,\pp)$.
Kim and Roush~(\cite[p.\ 92]{KR95}) proved:
\begin{prop}
  Let $k$ be a subfield of a $p$-adic field of odd residue
  characteristic $p$. Then $k$ has a finite extension $k'=k(i,a,\pp)$
  over which Hypothesis $\mathcal{H}$ is true.
\end{prop}

We will now show that for the purpose of proving
Theorem~\ref{Theorem}, we may enlarge the constant field $k$ (and
hence $\LL$).  In particular, we may assume that our field $k$
satisfies Hypothesis $\mathcal{H}$. Since we want to use the
coefficients of the quadratic form $q$ in our diophantine definitions,
we want to have $a$ and $\pp$ in our ring of coefficients.

\begin{prop}\label{coeffs2}
  Let $K$ be a field and let $A_0\subseteq K$ be a subring.\\
  (1) Let $c_1, \dots, c_m$ be elements of $K$ which are algebraic
  over $\Frac(A_0)$. If Hilbert's Tenth Problem for $K$ with
  coefficients in $A_0[c_1, \dots, c_m]$ is undecidable, then
  Hilbert's Tenth Problem for $K$ with coefficients in $A_0$ is
  undecidable.
  \\
  (2) Let $L/K$ be an extension which is generated by elements $b_1,
  \dots, b_r$ which are algebraic over $\Frac(A_0)$. If Hilbert's
  Tenth Problem for $L$ with coefficients in $A_0[b_1, \dots,b_r]$ is
  undecidable, then Hilbert's Tenth Problem for $K$ with coefficients
  in $A_0$ is undecidable.
\end{prop}
\begin{proof}
  (1) Since $c_i$ is algebraic over $\Frac(A_0)$ for $i = 1 , \dots,
  m$, we can consider its minimal polynomial over $\Frac(A_0)$. After
  multiplying by a common denominator we get an irreducible polynomial
  $p_i(x)$ over $\Frac(A_0)$ with coefficients in $A_0$. 
Given a polynomial equation $f(x_1, \dots, x_n)=0$ with coefficients
in $A_0[c_1, \dots, c_m]$, we can construct a system of polynomial equations
with coefficients in $A_0$ by replacing, for $i=1, \dots,m$, each
occurrence of $c_i$ in $f$ with a new variable $y_i$. Let $g(x_1,
\dots,x_n,y_1, \dots, y_m)$ be this new equation obtained from $f$, and
for $i=1 , \dots, m$, add the equation $p_i(y_i)=0$. Then the system of
equations $g=0, p_i(y_i)=0, i=1, \dots,m$, has a solution in $K$ if and
only if $f(x_1, \dots, x_n)=0$ has a solution in $K$. By
Lemma~\ref{combine} the system of equations can be replaced with one
single polynomial equation with coefficients in $A_0$.

(2) Since the $b_i$'s are algebraic over $\Frac(A_0)$, the minimal
polynomials of the $b_i$'s over $K$ have coefficients $d_1, \dots,
d_{\ell}$ which are algebraic over $\Frac(A_0)$. Now use the Technical
Lemma in Pheidas~\cite[p.\ 379]{Pheidas94} together with the first
part of this proposition.
\end{proof}
So in the following, whenever we pass to an extension $\LL' / \LL$,
we will choose the ring of coefficients $A_0$ large enough to ensure that
the elements generating $\LL' / \LL$ are algebraic over $A_0$.
We can work with an enlarged constant field $k$
that satisfies Hypothesis $\mathcal{H}$, since the elements $i,a,\pp$
specified there are algebraic over $\Q$.
\section{Proof of main theorem}\label{proof}
We need one more result from \cite{KR95} before we can prove
Theorem~\ref{Theorem}:
\begin{thm}\label{replacement}
  Let $k$ be a subfield of a $p$-adic field which satisfies hypothesis
  $(\mathcal{H})$, and let $a,\pp$ be as in hypothesis
  $(\mathcal{H})$. Let $g \in k(T)$ be such that $\ord_{T^{-1}}(g)=-2$.
  For $c_3,c_5 \in k$ let
  \[f(T)=f_{c_3,c_5}(T)= (1+T)^3g(T) + c_3T^3 + c_5T^5 .\] Let $U_0$
  be as in Theorem~\ref{densitytheorem}. If $\ord_T(g)=1$, then there
  exist $c_3,c_5 \in U_0$ such that the two quadratic forms
\begin{align}
&\langle T,Ta,-1,-f\rangle \langle 1,\pp\rangle \label{qf:1}\\
&\langle T,Ta,-1,-af\rangle \langle 1,\pp\rangle \label{qf:2}
\end{align}
are isotropic over $k(T)$ {\em(}and hence over any finite extension of
$k(T)${\em )}.
\end{thm}
\begin{proof}
 This theorem follows immediately
  from Theorem 9, Theorem 17, and Theorem 21 of \cite{KR95}.
\end{proof}  
\subsection{Proof for transcendence degree one}
We will now prove Theorem~\ref{Theorem} when $\LL$ is a finite
extension of the rational function field $k(\tau)$. Let $\kappa \subseteq k$ be
defined as in Section~\ref{kappa}.
\begin{thm}\label{special}
  Let $k$ be a subfield of a $p$-adic field of odd residue
  characteristic, and let $\LL$ be a finite extension of $k(\tau)$.
  There exists a finite set $\{c_1,\dots,c_{\ell}\}$ of elements of
  $\kappa(\tau)$, not all constant,
  such that Hilbert's Tenth Problem for $\LL$ with coefficients in
  $\Z[c_1,\dots,c_{\ell}]$ is undecidable.
\end{thm}
\begin{proof}
  Let $E_0,\EE_0,T$ be as in Section~\ref{ellipticsection}, {\it i.e.\ }
  $$\EE_0:(T^3+T+1)Y^2=X^3+X+1$$
  has rank one over $\LL$ with generator $P:=(T,1)$ (modulo 2-torsion)
  and there exists a prime $\QQ$ above $T^{-1}$ which is unramified.
  By Theorem~\ref{elliptic}, $T$ can be chosen to be algebraic
  over $\kappa(\tau)$.  After making a constant field extension as in
  Proposition~\ref{ff2} we may assume that the residue field of $\LL$
  at the prime $\QQ$ is $k$, and that $k$ is algebraically closed in
  $\LL$. After extending the constant field $k$ further, if necessary,
  we obtain an extension $k'$ that satisfies hypothesis $\mathcal{H}$.
  After these constant field extensions $\QQ$ remains unramified, and
  by Moret-Bailly's theorem (\cite[Theorem 1.8]{MB01}) the group
  $\EE_0(\LL)$ is still generated by $(T,1)$. Let $\LL':= \LL k'$. We
  apply Proposition~\ref{coeffs2}(2) to the finite extension $\LL'/
  \LL$, and we choose a ring of coefficients $A_0$ that satisfies the
  hypotheses of Proposition~\ref{coeffs2}, {\it i.e.}\ $A_0$ contains
  the coefficients of the minimal polynomial of $T$ over $\kappa(\tau)$ and
  the coefficients of the minimal polynomials of the elements
  generating $\LL' / \LL$.  By Proposition~\ref{ff2}, $A_0$ can be
  chosen to be of the form $A_0=\Z[c_1,\dots,c_{\ell}]$ with
  $c_1,\dots,c_{\ell} \in \kappa(\tau)$, with $\kappa$ as in Section~\ref{kappa}.
  Let $a, \pp$ be the elements of $k'$ as in Hypothesis $\mathcal{H}$
  and let $A:= A_0[T,a,\pp]$. The elements $a, \pp$ are algebraic over
  $\Q$. By Proposition~\ref{coeffs2}, proving
  that Hilbert's Tenth Problem for $\LL'$ with coefficients in $A$ is
  undecidable is enough to prove undecidability for $\LL$ with
  coefficients in $A_0$.  For simplicity of notation, we rename $\LL'$
  and $k'$ as $\LL$ and $k$ again, respectively.

  Let $P_m:=m(T,1)=(X_m,Y_m)$ and $\psi_m := X_m/TY_m$.  We will
  construct a diophantine model of $\langle \Z, 0,1 ;+,\cdot\rangle$
  in $\LL$ with coefficients in $A$.
  \\
  The elliptic curve $\EE_0$ is a projective variety, but any
  projective algebraic set can be partitioned into finitely many
  affine algebraic sets, which can then be embedded into a single
  affine algebraic set.  This implies that the set $\EE_0(\LL)$ is
  diophantine over $\LL$, since we can take care of the point at
  infinity $\mathbf{O}$ of $\EE_0$.  Hence the set
  \begin{align*}\mathcal{S'}&:=\{(X_{2n},Y_{2n}):n \in \Z\}\\&=\{(x,y)
\in \LL^2: \exists \,u,v \in \LL:(u,v) \in \EE_0(\LL) \land (x,y)=2(u,v)\}
\end{align*} is  diophantine over $\LL$ with coefficients in $\Z[T]$. 
Then the set 
 \begin{align*}\mathcal{S}:=\{&(X_{n},Y_{n}):n \in \Z\}\\=\{&(x,y)
   \in \LL^2: \exists \,n \in \Z:\\ &\left( (x,y)=(X_{2n},Y_{2n}) \lor
     (x,y)=(X_{2n},Y_{2n})+(T,1)\right)\}
\end{align*} is diophantine over $\LL$ with coefficients in $\Z[T]$ as well.

By associating the point $P_n=(X_{n},Y_{n})$ to an integer $n$ we obtain
a bijection between $\Z$ and $\mathcal{S}$, and addition of elements of
$\mathcal{S}$ is existentially definable, because it is given by the
group law on the elliptic curve.  It remains to show that
multiplication of elements of $\mathcal{S}$ is existentially
definable. Let $t:=T^{-1}$. We can consider $\LL$ as an extension of
$k(t)$. By the above discussion, the prime $\QQ$ above $t$ is
unramified, and the residue field of $\LL$ at $\QQ$ is $k$.

Let $q:=\langle 1,a\rangle \langle 1, \pp\rangle$ be the quadratic
form over $k$ as in hypothesis $\mathcal{H}$.

For $w \in \LL$ let $\Phi(w)$ be the formula expressing that the
quadratic forms
\begin{gather*}
  \langle t,-at,-1,-w\rangle \langle 1,\pp\rangle \text{ and } \langle
  t,-at,-1,-aw\rangle \langle 1, \pp\rangle
\end{gather*}
are isotropic over $\LL$. Clearly this is an existential formula. We
will show that $n \cdot m =r$ if and only of $\Phi(w)$ holds for a certain
function $w$ that is formed from the $x$- and $y$-coordinates
of the points $n\cdot (T,1)$, $m\cdot (T,1)$, and $r\cdot (T,1)$.

As before, given $n,m,r \in \Z-\{0,1,-1\}$ let $u:= \psi_{m}
\psi_{n}-\psi_{r} + (1/2) \cdot t$ and let $v:= \psi_{m} \psi_{n}-\psi_{r} + (1/3)
\cdot t$. Let $\ord_t,\ord_{t^{-1}}$ be the normalized
discrete valuations of $k(t)$ associated to $t$ and $t^{-1}$. By
Lemma~\ref{orders}, $\ord_{t^{-1}}(u)=-2$, $\ord_{t^{-1}}(v)=-2$ and
if $n\cdot m=r$, then $\ord_t(u) =1$ or $\ord_t(v)=1$. If $n\cdot m\neq r$,
then $\ord_t(u) =0$ and $\ord_t(v)=0$. (The cases where $n,m$ or $r$
are in $\{0,1,-1\}$ can be handled separately.)
Let $U_0$ be as in Theorem~\ref{densitytheorem}.
For $c_3,c_5 \in U_0$ let $$f_{(u,c_3,c_5)}:= (1+t)^3u + c_3t^3 +
c_5t^5,$$
and let 
$$f_{(v,c_3,c_5)}:= (1+t)^3v + c_3t^3 +c_5t^5.$$

We will show that 
\begin{equation}\label{mult}n \cdot m=r \leftrightarrow  \exists c_3,c_5 \in U_0: \left(
  \Phi(f_{(u,c_3,c_5)}) \lor \Phi(f_{v,c_3,c_5}) \right).\end{equation}
Since $U_0$ is diophantine over $\LL$ with coefficients in
$\Z$, it is easy to see that the condition that there exist $c_3,c_5
\in U_0$ for which the quadratic form $$\langle t,-at,-1,-f_{(u,c_3,c_5)}\rangle
\langle 1, \pp\rangle$$ has a solution in $\LL$ can be described by an
existential definition with coefficients in $A$.  Hence the
right-hand-side of (\ref{mult}) is
an existential definition with coefficients in $A$.

Suppose that $n\cdot m=r$. Then at least one of $u,v$ has order 1 at $t$.
Say $\ord_t(u)=1$. 
Then $f_{(u,c_3,c_5)}$ is an element of $k(t)$, and so
the quadratic forms
\begin{gather}\label{fnmr}
\langle t,-at,-1,-f_{(u,c_3,c_5)}\rangle \langle 1, \pp\rangle \text{ and }
\langle t,-at,-1,-af_{(u,c_3,c_5)}\rangle \langle 1, \pp\rangle
\end{gather}
are quadratic forms over $k(t)$.  By Theorem~\ref{replacement},
applied with $g=u$, there
exist $c_3,c_5 \in U_0$ such that $\Phi(f_{(u,c_3,c_5)})$ holds.

Conversely, assume that $n\cdot m\neq r$. Then by Lemma~\ref{orders}, $u$
and $v$ have order 0 at $t$, and so $f_{(u,c_3,c_5)}$ and
$f_{(v,c_3,c_5)}$ have order 0 at $t$ for any choice of $c_3,c_5 \in
U_0$.  Then by Lemma~\ref{anisotropic}, applied with
$g=f_{(u,c_3,c_5)}$, for any choice of $c_3,c_5 \in U_0$, one of the two
quadratic forms in (\ref{fnmr}) is anisotropic over $\LL$, so
$\Phi(f_{(u,c_3,c_5)})$ does not hold. Similarly, $\Phi(f_{(v,c_3,c_5)})$
does not hold for any choice of $c_3,c_5 \in U_0$.

\end{proof}
\subsection{Generalization to higher transcendence degree}
Let $k$ be a subfield of a $p$-adic field of odd residue
characteristic, and let $\LL$ be a finite extension of the rational
function field $k(\tau, \tau_2, \dots,\tau_n)$. Let $k_1$ be the algebraic
closure of $k(\tau_2, \dots, \tau_n)$ in $\LL$.  Then $\LL$ is a finite
extension of $k_1(\tau)$. Let $\kappa$ be as in Section~\ref{highertrdeg}.
We can apply Theorem~\ref{elliptic} to the elliptic curve $E_0$
defined in Section~\ref{ellipticsection} and with $\LL / k_1(\tau)$ to
obtain an element $T_1$ which is algebraic over $\kappa(\tau)$.  Consider
the elliptic curve $\EE_0$ defined by the affine equation
$(T_1^3+T_1+1)Y^2 =X^3 +X+1$. By Theorem~\ref{elliptic}, $\EE_0(\LL)$
is generated by $(T_1,1)$ (modulo 2-torsion).

To prove Theorem~\ref{Theorem} when the
transcendence degree of $\LL/k$ is $\geq 2$, we have to prove the
following Lemma.
\begin{lem}\label{anisotropic2}
  Let $k$ be a subfield of a $p$-adic field, and let $K$ be a
  finite extension of the rational function field $k(\tau_2,\dots,\tau_n)$.
  There exists a finite extension $k'/k$ such that $k'$ satisfies
  hypothesis $\mathcal{H}$, and such that the form $q=\langle1,a \rangle \langle 1,
  \pp \rangle$ as in Hypothesis $\mathcal{H}$ remains anisotropic over $K'
  =Kk'$.
\end{lem}
\begin{proof}
  Let $\PP_2$ be a prime of $K$ lying above the prime $\tau_2$ of the
  rational function field $k(\tau_3, \dots,\tau_n)(\tau_2)$, and let
  $k_{\PP_2}$ be the residue field of $\PP_2$.  Then $k_{\PP_2}$ is a
  finite extension of $k(\tau_3, \dots, \tau_n)$.  Now let $\PP_3$ be a
  prime of $k_{\PP_2}$ lying above the prime $\tau_3$ of $k(\tau_4,
  \dots,\tau_n)(\tau_3)$. Let $k_{\PP_3}$ be the residue field of $\PP_3$.
  The field $k_{\PP_3}$ is a finite extension of $k(\tau_4,\dots,\tau_n)$.
  After repeating this process we obtain a finite extension
  $k_{\PP_n}$ of $k$.  
  From the proof of \cite[Proposition 8]{KR95} it follows that we can
  find a finite extension $k'$ of $k$ which is generated by elements
  algebraic over $\Q$ such that both $k'$ and $k_{\PP_n}k'$ satisfy
  Hypothesis $\mathcal{H}$.
  \\
  {\bf Claim:} The field $k'$ has the desired property.\\
  {\bf Proof of Claim:} Let $\PP_2'$ be a prime of $K'$ extending
  $\PP_2$. Let $k_{\PP_2'}$ be the residue field of $\PP_2'$. Then
  $k_{\PP_2'}= k'k_{\PP_2}$ by Theorem~\ref{resextension}. Let
  $\PP_3'$ be a prime of $k_{\PP_2'}$ extending $\PP_3$, and let
  $k_{\PP_3'}$ be the residue field.  Define $\PP_4', k_{\PP_4'},
  \dots, \PP_n',k_{\PP_n'}$ similarly. We have $k_{\PP_n'} =
  k'k_{\PP_n}$.  Assume by contradiction that $q$ is isotropic over
  $K'$. Take a solution $f_1,\dots,f_4$ in $K'$.  By scaling the
  $f_i$'s with the same factor we can arrange it so
  that $\ord_{\PP_2'}f_i \geq 0$ for $i=1, \dots, 4$, and by changing
  the $f_i$'s further we may assume that $\ord_{\PP_2'}f_i =0$ for
  some $i$. Now look at $\overline{f_i}=f_i \mod \PP_2'$, $i=1, \dots,
  4$. This gives us that $q$ is isotropic over $k_{\PP_2'}$. By
  repeating this we obtain that $q$ is isotropic over
  $k_{\PP_n'}$, contradiction.
\end{proof}
Now we can generalize the proof of Theorem~\ref{special} and prove
Theorem~\ref{Theorem}:

\begin{proof}[{\bf Proof of Theorem~\ref{Theorem}}]
   Let $k, k_1, \LL, E_0,\EE_0$ be as above. The proof proceeds as in
   Theorem~\ref{special}. Let $\QQ$ be a prime of $\LL$ above the
   prime $T_1^{-1}$ of $ k(\tau_2,\dots,\tau_n)(T_1)$ which is unramified.
   Let $m(T_1,1)=(X_m,Y_m)$ and $\psi_m := X_m/T_1 Y_m$.  As in
   Remark~\ref{higher}, we enlarge $k_1$ to a finite extension $k_1'$
   such that the prime $\QQ'$ of the compositum $\LL k_1'$ (extending
   $\QQ$) has residue field $k_1'$.  Let $\LL':= \LL k_1'$.
   We now pass to an extension $k'$ of the constant field $k$ (and
   hence enlarge $k_1'$ and $\LL'$ further) such that $k'$ is as in
   Lemma~\ref{anisotropic2} for the extension $k_1' / k(\tau_2, \dots,
   \tau_n)$.  By the proof of Lemma~\ref{anisotropic2}, $k'/k$ can be
   generated by elements which are algebraic over $\Q$.  By
   \cite[Proposition 8.3]{Rosen}, $k_1'$ is algebraically closed in
   $\LL'$.

   Now we choose our ring of coefficients $A_0$ such that
   $\LL' / \LL$ is generated by elements algebraic over $A_0$, and
   such that $A_0$ contains the coefficients of the minimal polynomial
   of $T_1$.
   By the above arguments, together with Remark~\ref{higher} and
   Theorem~\ref{elliptic}, $A_0$ is of the form $A_0= \Z[
   c_1, \dots, c_{\ell}]$, with $\{c_1,\dots,c_{\ell}\} \in \kappa(\tau)$. We
   let $A:=A_0[T_1,a,\pp]$, with $a, \pp$ the elements as in
   Hypothesis~$\mathcal{H}$. For simplicity of notation we rename
   $\LL',k_1',k'$, and $\QQ'$ as $\LL, k_1,k$, and $\QQ$.

   By the remark after Theorem~\ref{densitytheorem} we can still
   construct $U_0 \subseteq k$ which is diophantine over $\LL$ with
   coefficients in $\Z$ and whose intersection with $\Q$ is dense in
   any finite product of $p$-adic fields.

   To prove the theorem we will construct a diophantine model of the
   structure $\langle \Z,0,1;+,\cdot\rangle$ in $\LL$ with
   coefficients in $A$.  As before let $\mathcal{S}:=\{(X_{n},Y_{n}):n
   \in \Z\}$. This set is diophantine over $\LL$ with coefficients in
   $\Z[T_1]$, and it remains to show that multiplication of elements
   of $\mathcal{S}$ is existentially definable. Let $t_1:=T_1
   ^{-1}$. Then $\LL$ is a finite extension of $k_1(t_1)$ and the
   prime $\QQ$ of $\LL$ above $t_1$ is unramified.

Let $q:=\langle 1,a\rangle \langle 1, \pp\rangle$ be the quadratic
form over $k$ as in hypothesis $\mathcal{H}$.

For $w \in \LL$ let $\Phi(w)$ be the formula expressing that the
quadratic forms
\begin{gather}\label{fnmr2}
  \langle t_1,-at_1,-1,-w\rangle \langle 1,\pp\rangle \text{ and } \langle
  t_1,-at_1,-1,-aw\rangle \langle 1, \pp\rangle
\end{gather}
are isotropic over $\LL$.
Given $n,m,r \in \Z-\{0,1,-1\}$ let $u:= \psi_{m} \psi_{n}-\psi_{r} + (1/2) \cdot t_1$
and $v:= \psi_{m} \psi_{n}-\psi_{r} + (1/3) \cdot t_1$.  For $c_3,c_5 \in k$ let
$$f_{(u,c_3,c_5)} := (1+t_1)^3u + c_3t_1^3 + c_5t_1^5$$ and
\[
f_{(v,c_3,c_5)} := (1+t_1)^3v + c_3t_1^3 +
c_5t_1^5
\]
The elements $f_{(u,c_3,c_5)}$ and $f_{(v,c_3,c_5)}$
are elements of $k(t_1)$.\\
 
We will show that 
\begin{equation}\label{eq:5}n \cdot m =r \leftrightarrow \exists c_3,c_5 \in U_0:
\left( \Phi(f_{(u,c_3,c_5)} \lor \Phi(f_{(v,c_3,c_5)})\right).\end{equation}
The same argument as in Theorem~\ref{special} shows that the
right-hand-side of (\ref{eq:5}) is existential with coefficients in $A$.
Suppose that $n\cdot m=r$. Then at least one of $u, v$ has order 1 at
$t_1$. Say $\ord_{t_1}(u)=1$. Then $f_{(u,c_3,c_5)}$ is an element of
$k(t_1)$, and the same argument as in Theorem~\ref{special} shows that
there exist $c_3,c_5 \in U_0$ such that $\Phi(f_{(u,c_3,c_5)})$ holds.

Conversely, assume that $n \cdot m\neq r$. Then by Lemma~\ref{orders}, $u$
and $v$ have order 0 at $t_1$, and so
$f_{(u,c_3,c_5)},f_{(v,c_3,c_5)}$ have order 0 at $t_1$ for any choice
of $c_3,c_5 \in U_0$.  Then by Lemma~\ref{anisotropic2}, we can apply
Lemma~\ref{anisotropic} to the extension $\LL / k_1(t_1)$, and with
$g=f_{(u,c_3,c_5)}$. Hence for any choice of $c_3,c_5 \in U_0$,
$\Phi(f_{(u,c_3,c_5)})$ does not hold. Similarly, for any choice of
$c_3,c_5 \in U_0$,  $\Phi(f_{(v,c_3,c_5)})$ does not hold, either.
\end{proof}
\section{Appendix}
In this section we will state the definitions and theorems about
quaternion algebras and quadratic forms that we used in our proof.
We need the following two definitions.
\begin{defn}\label{quaternion}
  Let $F$ be a field of characteristic $\neq 2$, and let $a,b \in
  F^{*}$. We define the {\em quaternion algebra} $\left(
    \frac{a,b}{F}\right)$ to be the $F$-algebra on two generators
  $i,j$ with defining relations: $i^2 =a$ and $j^2 =b$, and $ij=-ji$.
  The {\em associated norm form} of the quaternion algebra $\left(
    \frac{a,b}{F}\right)$ is the quadratic form $\langle
  1,-a,-b,ab\rangle$.
\end{defn}
\begin{defn}\label{def2}
  We say that a quaternion algebra $A=\left( \frac{a,b}{F}\right)$
  {\it splits over $F$} if $A \isom M_2(F)$.
\end{defn}
We can see whether a quaternion algebra is split by looking at its
norm form:
\begin{prop}\label{split}
The quaternion algebra $\left( \frac{a,b}{F}\right)$ splits over $F$
iff its associated norm form is isotropic.
\end{prop}
\begin{proof}
This is proved in \cite[Theorem 2.7, p.\ 58.]{Lam}
\end{proof}
\begin{prop}\label{multiply}
For $a,b,c \in F^{*}$, we have
\[
\left( \frac{a,b}{F}\right) \tensor \left(\frac{a,c}{F} \right)\isom
\left( \frac{a,bc}{F}\right) \tensor M_2(F).
\]
\end{prop}
\begin{proof}
This is Corollary 2.11 in \cite[p.\ 61]{Lam}.
\end{proof}

{\bf{Acknowledgments.}} I thank Brian Conrad and Thanases Pheidas for
several helpful discussions. I also thank Laurent Moret-Bailly for
suggesting a way to simplify the proof of Theorem~\ref{help}.

\end{document}